\documentclass[twoside]{article}
 \usepackage{graphicx,amsmath,latexsym,amssymb,amsthm,fancyhdr}

\font\chuto=cmbx10 at 14pt
\font\chudaude=cmcsc10

\font\ita=cmti9
 \numberwithin{equation}{section}

 \newtheorem{thm}{Theorem}[section]
\newtheorem{prop}[thm]{Proposition}
\newtheorem{lem}[thm]{Lemma}
\newtheorem{cor}[thm]{Corollary}
\newtheorem{rem}[thm]{Remark}
\newtheorem{definition}[thm]{Definition}
 
\def\beginpf{\noindent{\bf Proof.\;}}
\def\endpf{\hfill$\Box$}
\title{ 
{\chuto   TOWARDS QUANTITATIVE CLASSIFICATION 
          OF CAYLEY AUTOMATIC GROUPS}
          \footnotetext{ 
 \\
  \\
  {\bf Key
words:} automatic group, Cayley automatic group, nilpotent group, fundamental group, 
Dehn function, growth function, numerical characteristic.
 }
\author{{\bf Dmitry Berdinsky${}^{*}$ {\rm and} Phongpitak Trakuldit${}^{\dagger}$    
}\\
 \vspace{-0.2cm}\\
 \vspace{-0.2cm} {\ita   Department of   
 Mathematics, Faculty of Science, Mahidol University
} \\
 \vspace{-0.2cm} {\ita  Centre of Excellence in Mathematics, Commission on Higher Education}\\
\vspace{-0.2cm}{\ita Bangkok, 10400, Thailand  
}\\
 \vspace{-0.2cm} {\ita email: ${}^{*}$berdinsky@gmail.com,${}^{\dagger}$p.trakuldit@gmail.com }}}

\begin{document}
\label{firstpage}

\date{ }

\maketitle

\begin{abstract}
 In this paper we 
 address the problem of quantitative 
 classification of Cayley automatic groups 
 in terms of a certain numerical characteristic 
 which we earlier introduced for this  
 class of groups. 
 For this numerical characteristic
 we formulate and prove a fellow traveler property, 
 show its relationship with the Dehn function 
 and prove its invariance with respect to
 taking finite extension, direct product and 
 free product. 
 We study this characteristic for 
 nilpotent groups with a particular
 accent on the Heisenberg group, 
 the fundamental groups of torus bundles over the
 circle and groups of exponential growth. 
\end{abstract}

\section{\bf Introduction and Preliminaries}

   Strings over a finite alphabet appear 
   a natural way to represent elements 
   of a finitely generated group. 
   Following this way 
   Thurston introduced automatic 
   groups which became an 
   important part of geometric 
   group theory \cite{Epsteinbook}. 
   Trying to extend
   the class of automatic groups, 
   one can either use more powerful 
   computational models (e.g., asynchronous 
   automata, pushdown automata and etc.) 
   or relax the constraint on 
   the correspondence between 
   strings and group elements 
   (for automatic groups this 
   correspondence is given by the 
   canonical map).   
   The latter approach leads to 
   Cayley automatic groups introduced 
   by Kharlampovich, Khoussainov and 
   Miasnikov~\cite{KKM11}. 
   Utilization of both approaches   
   simultaneously 
   leads further  to 
   $\mathcal{C}$--graph automatic groups 
   introduced by Elder and Taback
   \cite{ElderTabackCgraph}. In 
   this paper we focus only on Cayley 
   automatic groups. 
   
   Cayley automatic groups utilize 
   exactly the same computational 
   model as automatic groups, 
   so they preserve some key 
   algorithmic features of 
   automatic groups, but
   the correspondence between 
   strings and group elements 
   can be arbitrary.  
   Another way to define Cayley automatic 
   groups is to say that they 
   are finitely generated groups 
   for which labeled directed Cayley 
   graphs are automatic 
   (FA--presentable) structures 
   \cite{KhoussainovNerode95,KM07,KN08}. 
   For a recent survey 
   of the theory of automatic 
   structures we refer the reader 
   to~\cite{Stephan2015}. 
   The class of Cayley automatic groups 
   is essentially wider than 
   the class of automatic groups 
   \cite{KKM11}.    
   Also, Cayley automatic groups include 
   important classes of groups 
   such as nilpotent groups of
   nilpotency class two, 
   fundamental groups of $3$--manifolds, 
   Baumslag--Solitar groups, 
   restricted wreath products  
   of Cayley automatic groups by 
   the infinite cyclic group, 
   higher rank lamplighter groups  
   \cite{KKM11,dlt14,Taback18}. 
  
   We assume that the reader is familiar 
   with the definitions of finite 
   automata and regular languages (a concise 
   introduction is given in, e.g., 
   \cite[Sections~1.1--2]{Epsteinbook}). 
   For a given finite alphabet $\Sigma$ we 
   denote by $\Sigma^*$ the set of all 
   finite strings over $\Sigma$ and by 
   $\Sigma_\diamond$ 
   the alphabet 
   $\Sigma = \Sigma \cup \{ \diamond \}$ 
   (it is assumed that $\diamond \not\in \Sigma$). 
   For any $w \in \Sigma^*$, we denote by 
   $|w|$ the length of the string $w$.  
   Let $w_1, \dots ,w_n 
        \in \Sigma^*$.
   The convolution $w_1 \otimes \dots  \otimes 
   w_n$ is the string of a length 
   $m = \max \{|w_1|, \dots, |w_n|\}$ over the 
   alphabet ${\Sigma_\diamond ^n}' =  
   \Sigma_\diamond ^n \setminus 
   \{(\diamond,\dots,\diamond)\}$ for which 
   the $k$th symbol, $k=1,\dots,m$, is 
   $(\sigma_{1k},\dots,\sigma_{nk}) \in 
   {\Sigma_\diamond ^n}' $, 
   where $\sigma_{ik}$ 
   is the $k$th symbol 
   of $w_i$ if $k \leqslant |w_i|$ and 
   $\sigma_{ik}= \diamond$ if $k>|w_i|$ 
   for $i=1,\dots,n$.  
   For any  relation 
   $R \subseteq \Sigma^{*n}$, 
   we say that $R$ is FA--recognizable (regular) 
   if $\otimes R = \{ w_1 \otimes \dots 
   \otimes w_n \, | \, 
   (w_1,\dots,w_n) \in R \}$ is a regular 
   language over the alphabet 
   ${\Sigma_\diamond ^n}'$.                  
     Let $G$ be a finitely generated (f.g.) 
     group and $A \subset G$ be a finite 
     generating set of $G$. 
     Let $A^{-1}$ be the set of 
     the inverses of elements of $A$ 
     and $S = A \cup A^{-1}$.       
     We denote by $\pi: S^* \rightarrow G$ 
     the canonical map 
     which maps any given 
     string $w=s_1 \dots s_n \in S^*$  
     to the group element 
     $g=s_1 \dots s_n \in G$.     
   \begin{definition} 
   \label{cayleyautomatic_def1}      
       A group $G$ is called Cayley automatic 
       if there exists a bijection 
       $\psi: L \rightarrow G$       
       between 
       some regular language 
       $L \subseteq \Sigma^*$ and the group $G$ for 
       which the binary relation 
       $R_a = \left\{ \left(\psi^{-1}(g),
        \psi^{-1}(ga)
       \right)
        | g \in G \right\}$ is FA--recognizable for 
        every $a \in A$.    
        Such a bijection 
        $\psi: L \rightarrow G$ is called 
        a Cayley automatic representation 
        of $G$.  
    \end{definition}           
       In this paper we assume that $\Sigma = S$,  
       unless otherwise stated. 
       This assumption is needed 
       to correctly define the function 
       $h(n)$ in the formula \eqref{def_func_h_1} 
       below:
       if $w \in S^*$, then $\pi(w)$ is 
       in the group $G$ as well as $\psi(w)$, 
       so one can get the 
       distance 
       $d_A(\pi(w),\psi(w))$       
       between $\pi(w)$ and $\psi(w)$ 
       in the Cayley graph $\Gamma (G,A)$.                 
       We recall that for given 
       $g_1,g_2 \in G$, the distance 
       $d_A (g_1,g_2)$       
       between the elements
        $g_1$ and $g_2$ in $G$
       with respect to $A$ is the length of a
       shortest path from $g_1$ to 
       $g_2$ in the Cayley graph 
       $\Gamma(G,A)$. 
       For a given $g \in G$, 
       we denote by 
       $d_A (g)$ the distance 
       $d_A (e,g)$, where $e$ 
       is the identity of the 
       group $G$.  
       Since the cardinality  
       of $S$ is at least two, 
       it can be verified 
       that Definition 
       \ref{cayleyautomatic_def1} 
       (either together with the assumption 
       that $\Sigma = S$ or without it)       
       is equivalent to the original 
       definition of Cayley automatic groups 
       \cite[Definition~6.4]{KKM11}       
       (they are also referred as
        Cayley graph automatic 
        or graph automatic groups 
        in the literature).
       Furthermore, assuming that 
       $\Sigma = S$ and $\psi = \pi$ in 
       Definition \ref{cayleyautomatic_def1}, 
       one gets the definition of automatic 
       groups; it can be also verified that 
       it is equivalent to the original 
       definition given by Thurston, see
       \cite[Definition~2.3.1]{Epsteinbook}.                            
       This observation motivated us to 
       introduce a function
       \eqref{def_func_h_1} 
       as a measure of deviation of 
       a given Cayley automatic representation 
       $\psi$ from automatic representations
       \cite{measuring_closeness1}:
       \begin{equation} 
       \label{def_func_h_1} 
       h(n) = \max \left\{d_A (\pi(w),\psi(w))| 
       w \in L ^{\leqslant n}\right\}, 
       \end{equation}
       where 
       $L^{\leqslant n} = \{ 
       w \in L\,|\, |w| \leqslant n \}$
       is the set of strings from 
       $L$ of a length less or equal than $n$.     
       If a group $G$ is Cayley automatic 
       but not automatic, a Cayley automatic 
       representation $\psi$ for which 
       $\psi = \pi$ does not exist. So, in 
       this case, 
       for every Cayley 
       automatic representation $\psi$ of $G$
       the function $h(n)$ defined by 
       \eqref{def_func_h_1} 
       is  not identically 
       equal to zero.

       We  denote by $\mathfrak{F}$ the
       set of all nondecreasing functions 
       from some interval 
       $[Q,+\infty) \subseteq \mathbb{N}$ to 
       the set of nonnegative real numbers.         
       Clearly, a function $h(n)$ given in 
       \eqref{def_func_h_1} is in 
       $\mathfrak{F}$. For any given 
       $g,f \in \mathfrak{F}$, we say that  
       $g \preceq f$  ($g$ is coarsely less or
       equal than $f$) if there exist
       nonnegative integer $N$ and 
       positive integers $K$ and $M$ 
       for which $g(n) \leqslant K f(Mn)$ for 
       all $n \geqslant N$. We say that
       $g \asymp f$ ($g$ is coarsely equal
       to $f$) if $g \preceq f$ and 
       $f \preceq g$. Similarly, 
       we say that $g \prec f$ 
       ($g$ is coarsely strictly less than $f$) 
       if $g \preceq f$ and $g \not\asymp f$.
       Clearly, the coarse equality 
       $\asymp$ gives an equivalence 
       relation on $\mathfrak{F}$. In this 
       paper we will be considering  
       functions from $\mathfrak{F}$ up to 
       this equivalence relation.

       Any given Cayley automatic group $G$ 
       admits infinitely many Cayley 
       automatic representations 
       $\psi: L \rightarrow G$.  
       So, in general, the problem of finding 
       Cayley automatic representations  
       minimizing coarsely the function 
       \eqref{def_func_h_1} is nontrivial. 
       In \cite[Theorems~11 and 13]
       {measuring_closeness1}, 
       we constructed Cayley automatic 
       representations of the 
       Baumslag--Solitar groups 
       $BS(p,q)$, $q > p \geqslant 1$ and 
       the lamplighter group 
       $\mathbb{Z}_2 \wr \mathbb{Z}$ 
       which are  minimizers of 
       the function \eqref{def_func_h_1}. 
       In both cases the minimum
       for the function $h(n)$ 
       is the identity function 
       $\mathfrak{i}$: 
       $\mathfrak{i} (n) = n$ for all 
       $n \in \mathbb{N}$.
       Furthermore, in \cite{measuring_closeness1} 
       we introduced classes of Cayley 
       automatic groups 
       $\mathcal{B}_f$ as follows. 
       For a given $f \in \mathfrak{F}$, 
       $G \in \mathcal{B}_f$ 
       if there exists 
       a Cayley automatic representation
       $\psi: L \rightarrow G$ 
       for which $h \preceq f$, where 
       $h$ is given by \eqref{def_func_h_1}. 
       In particular, the 
       Baumslag--Solitar groups 
       $BS(p,q)$, $q > p \geqslant 1$ and 
       the lamplighter group 
       $\mathbb{Z}_2 \wr \mathbb{Z}$ 
       are in the class $\mathcal{B}_\mathfrak{i}$
       and they cannot be in any 
       class $\mathcal{B}_f$ 
       if $f \prec \mathfrak{i}$.     
       
       It is easy to show that the definition 
       of a class $\mathcal{B}_f$ 
       does not  depend on the choice of generators
       \cite[Proposition~5]{measuring_closeness1}.
       Clearly, 
       $\mathcal{B}_f \subseteq 
       \mathcal{B}_g$ if $f \preceq g$.        
       Also, for the zero function $\bf{z}$, 
       where ${\bf z}(n)=0$ for all 
       $n \in \mathbb{N}$, the class 
       $\mathcal{B}_{\bf z}$ coincides with 
       the class of automatic groups. 
       In
       \cite[Theorem~8]{measuring_closeness1}
       we proved that 
       there exists no nonautomatic 
       group in any class $\mathcal{B}_d$, 
       where $d \in \mathfrak{F}$ is a 
       function bounded from above by 
       some constant; that is, 
       $\mathcal{B}_d = \mathcal{B}_{\bf z}$
        for
       any such function $d$. 
       Another group that we considered in 
       \cite{measuring_closeness1} was 
       the Heisenberg group 
       $\mathcal{H}_3 (\mathbb{Z})$.  
       We showed that $\mathcal{H}_3 
       (\mathbb{Z}) \in \mathcal{B}_\mathfrak{e}$,
       where $\mathfrak{e}$ is the exponential 
       function: $\mathfrak{e}(n)= \exp(n)$. 
       But a lower bound for $h(n)$ which 
       we could find in the case of 
       $\mathcal{H}_3 (\mathbb{Z})$  is 
       far from 
       being exponential, it is $\sqrt[3]{n}$
       \cite[Theorem~15]{measuring_closeness1}.
       
       For a given $G \in \mathcal{B}_f$ 
       we treat  
       $f \in \mathfrak{F}$ as a 
       numerical characteristic 
       of $G$.  We especially interested in 
       those $f$ which are sharp lower 
       bounds for \eqref{def_func_h_1}.
       The fact that the sharp lower 
       bounds
       can be obtained for some groups 
       sounds promising.        
       Numerical characteristics 
       of groups, e.g. growth functions, 
       Dehn functions, drifts of simple 
       random walks and etc., and 
       relations between 
       them are very important in group 
       theory, see, e.g., \cite{Vershik99}. 
       Another motivation to study 
       this numerical characteristic 
       is to address the problem 
       of characterization of Cayley 
       automatic groups; see also 
       \cite{dlt16}, where this 
       problem is addressed in terms of  
       numerical characteristics of 
       Turing transducers.    
             
       In this paper we continue 
       studying this numerical 
       characteristic of Cayley 
       automatic groups and its 
       relation to other numerical  
       characteristics initiated in  
       \cite{measuring_closeness1}. 
       In Section \ref{fellowdehnsec} 
       we propose a fellow 
       traveler property for 
       Cayley automatic groups in 
       Theorem \ref{fellow_theorem1}
       and show a relation with the Dehn 
       function in Theorem   
       \ref{isoperimetric_theorem1}.     
       The fellow traveler 
       property is well known for automatic
       groups but 
       its analog for Cayley automatic 
       groups had not been formulated before.
       In Section \ref{extfreeproddirprod} 
       we prove invariance of classes 
       $\mathcal{B}_f$ under taking finite 
       extension, direct product and 
       free product in Theorems 
       \ref{finiteindex_theorem1}, 
       \ref{directprodthm1} and 
       \ref{freeprodtheorem1}, respectively; 
       in the latter case we require 
       the function $f$ to satisfy a certain 
       inequality. 
       
       In Section \ref{nilpotent_section} 
       we show that the semidirect
       products 
       $\mathbb{Z}^n \rtimes_A \mathbb{Z}$, 
       unitriangular matrix groups 
       $UT_n (\mathbb{Z})$ and  
       all f.g. nilpotent groups of 
       nilpotency class two are in the 
       class $\mathcal{B}_\mathfrak{e}$, 
       see Theorem \ref{nilpotentsec_thm1}. 
       However, 
       this result       
       is obtained from 
       certain Cayley automatic 
       representations of these 
       groups and we do not know whether 
       they are minimizers of 
       the function  \eqref{def_func_h_1} 
       or not. We partly address this 
       issue in Theorem 
       \ref{nilpotentsec_starredlangthm1}
       by showing that if  a virtually 
       nilpotent group $G$ is in a class 
       $\mathcal{B}_p$ for some polynomial 
       $p$, then the language $L$ 
       of a Cayley automatic representation 
       $\psi: L \rightarrow G$, for which 
       $h \preceq p$, must be simply starred. 
       
       In Section \ref{nonautomaticitysection} 
       we address the problem of sharp 
       lower bounds of the function
       \eqref{def_func_h_1} specifically for the 
       Heisenberg group 
       $\mathcal{H}_3 (\mathbb{Z})$. 
       In Theorem \ref{nonautth1} we show 
       that under a certain condition  
       on a Cayley automatic representation
       of $\mathcal{H}_3 (\mathbb{Z})$ 
       the growth of the function  
       \eqref{def_func_h_1} must be 
       at least exponential. 
       We note that the proof of 
       Theorem \ref{nonautth1} does not
       use any knowledge about 
       growth of the Dehn function, 
       which is very often used to 
       show that a given group 
       is not automatic. 
       We believe that 
       Theorem \ref{nonautth1} can 
       be useful for proposing 
       new approaches to proving 
       nonautomaticity of groups. 
       Section \ref{linearupperalmostall} 
       concludes the paper by showing 
       that for any Cayley 
       automatic representation $\psi: 
       L \rightarrow G$ of a group 
       of exponential growth    
       a linear upper bound 
       $d_A (\pi (w), \psi(w)) \leqslant 
       C |w|$ holds for almost all $w \in L$ 
       in a certain sense, see Theorem 
       \ref{allmostall_quasigeodesic_theorem1}. 
       However, in  Remark 
       \ref{exptowergrowthremark} 
       we explain that one should 
       be careful with 
       this simple observation
       made in Theorem  
       \ref{allmostall_quasigeodesic_theorem1}  by 
       constructing Cayley 
       automatic representations of 
       the lamplighter group 
       $\mathbb{Z}_2 \wr \mathbb{Z}$ for which
       the function \eqref{def_func_h_1}    
       grows faster than any tower of 
       exponents.  
       
       All questions that we 
       posed in 
       \cite[\S7]{measuring_closeness1}
       remain open. Let us pose 
       an additional question here: 
       is there any Cayley automatic 
       representation of a group of 
       polynomial growth 
       (which is not virtually abelian) 
       or a fundamental group 
       of a $3$--manifold (which is not
       automatic) for which the function 
       \eqref{def_func_h_1} is coarsely strictly  
       less than the exponential function 
       $\mathfrak{e}$?                   
       
 \section{\bf Fellow Traveler Property 
          and Connection with Dehn Functions}   
 \label{fellowdehnsec}
 In this section  we formulate a fellow traveler 
 property for Cayley automatic groups and 
 obtain a relation between 
 the Dehn function of a group 
 $G \in \mathcal{B}_f$ and 
 a function $f$.   
 For any word $w \in S^*$ and 
 nonnegative integer $t$ we put 
 $w(t)$ to be the prefix of $w$ of a
 length $t$ if $t\leqslant |w|$ and 
 $w$ if $t>|w|$. We denote 
 by $\widehat{w}: 
 [0,\infty) \rightarrow \Gamma (G,A)$ 
 the corresponding 
 path in the Cayley graph $\Gamma (G,A)$: 
 if $t$ is an integer, 
 then $\widehat{w}(t)= \pi(w(t))$ and
 if $t$ is not an integer, 
 $\widehat{w}(t)$ is obtained 
 by moving along the edge 
 $\left(\widehat{w}(\lfloor t \rfloor),
  \widehat{w}(\lceil t \rceil)\right)$ 
  with unit speed; 
  we will use only integer 
 values of $t$.  
 Let $\psi: L \rightarrow G$ be 
 any Cayley automatic representation 
 of a group $G$. We denote by $s$ be the following function:
 \begin{equation}
 \label{def_function_f_1}  
  s(n)= \max \{d_A (\widehat{w_1}(t),
  \widehat{w_2}(t))| \psi(w_1)g = \psi(w_2),
  g \in A, t \leqslant n\}.
 \end{equation} 
 That is, for every two words 
 $w_1,w_2 \in L$ representing 
 neighboring vertices in the 
 Cayley graph  $\Gamma(G,A)$ 
 (i.e., for some
  $g \in A$, $\psi(w_1)g=\psi(w_2)$) 
  the distance between 
 $\widehat{w_1}(t)$ and $\widehat{w_2}(t)$ 
 for all $t \leqslant n$ 
 is bounded from above by $s(n)$.  If $G$ is automatic 
 and $\psi$ is an automatic 
 representation of $G$, then $s(n)$ must be 
 a bounded function due to the 
 fellow traveler property for automatic 
 groups \cite[Lemma~2.3.2]{Epsteinbook}.

 \begin{thm} 
 \label{fellow_theorem1}        
    Assume that $G \in \mathcal{B}_f$ for 
    some nonzero function $f \in \mathfrak{F}$.  
    Then there is a Cayley automatic 
    representation 
    $\psi: L \rightarrow G$ such that
    for the function $s(n)$ given by 
    \eqref{def_function_f_1}, 
    $s \preceq f$.    
 \end{thm} 
 \beginpf 
   Since $G \in \mathcal{B}_f$, 
   there exists a Cayley automatic 
   representation 
   $\psi: L \rightarrow G$ such that 
   for the function 
   $h(n) = \max \{ d_A (\pi(w),\psi(w)) | 
    w \in L^{\leqslant n} \}$, $h \preceq f$.  
   Let $t,n$ be some nonnegative 
   integers for which $t \leqslant n$ and 
   $w_1,w_2 \in L$ be some words representing
   neighboring vertices in $\Gamma(G,A)$ 
   (i.e., $\psi(w_1)g =\psi (w_2)$ for some 
   $g \in A$). 
   The convolution $w_1 \otimes w_2$ 
   is in a regular language 
   $\otimes R_g$ accepted by 
   some two--tape synchronous automaton $M_g$.
   Let $T$ be a maximal 
   number of states in the automata 
   $M_g$ for all $g \in A$.  
   We assume that $t>T$. 
   If $t \leqslant \max\{|w_1|,|w_2|\}$, 
   there exist strings $u_1,v_1,u_2,v_2$ 
   for which 
   $u_1$ and $u_2$ are prefixes 
   of $w_1(t)$ and 
   $w_2(t)$ such that 
   $|u_1|,|u_2| \geqslant t - T$ and  
   for the strings 
   $w_1'=u_1 v_1$ and $w_2'=u_2 v_2$, 
   $|w_1'|,|w_2'|\leqslant t$   
   and 
   the convolution $w_1' \otimes w_2' \in
   \otimes R_g$.
   If $t \geqslant \max \{|w_1|,|w_2|\}$, 
   then we simply put $u_1=w_1$, 
   $u_2=w_2$ and $v_1=v_2=\epsilon$, 
   where $\epsilon$ is the empty 
   string.    
   We have: $d_A (\widehat{w_1}(t),
                 \widehat{w_2}(t)) 
                 \leqslant
            d_A (\pi(u_1),\pi(u_2))+ 
            2T  \leqslant 
            d_A (\pi(w_1'),\pi(w_2')) + 
            |v_1| + |v_2| + 2T  \leqslant
            d_A (\pi(w_1'),\pi(w_2')) + 4T              
            $.  
    Moreover, 
    $d_A(\pi(w_1'),\pi(w_2'))
     \leqslant d_A (\pi(w_1'),\psi(w_1')) 
     + d_A (\psi(w_1'),\psi(w_2')) + 
       d_A (\psi(w_2'),\pi(w_2'))
     \leqslant  
     h(|w_1'|) + 1 + h(|w_2'|) 
     \leqslant 2 h(t) + 1         
    $. 
    Therefore, $d_A (\widehat{w_1}(t),
                 \widehat{w_2}(t)) 
                \leqslant 
                2 h(t) + 4T + 1 \leqslant 
                2 h (n) + 4T + 1$. 
    If $t \leqslant T$, 
    $d_A(\widehat{w_1}(t),\widehat{w_2}(t))$
    can be bounded from above by 
    $2T$.      
    Since $h \preceq f$ 
    and $f$ is a nonzero function, 
    then $s \preceq f$.
 \endpf 
 \begin{rem}
    Clearly, we have $d_A(\widehat{w_1}(t), 
    \widehat{w_2}(t)) \leqslant 
    d_A (\widehat{w_1}(t)) + 
    d_A (\widehat{w_2}(t))\leqslant 2t$. 
    Therefore, $s \preceq \mathfrak{i}$ 
    for any function $s$ given by 
    \eqref{def_function_f_1}. 
    So, Theorem \ref{fellow_theorem1} is  
    of interest if $f \prec \mathfrak{i}$.   
    It is not known whether there exists
    any Cayley automatic group in a class 
    $\mathcal{B}_f$, for 
    $f \prec \mathfrak{i}$, which 
    is not automatic.   
    If such groups do not exist, 
    Theorem \ref{fellow_theorem1} 
    might be a first step to prove it. 
    At least, Theorem \ref{fellow_theorem1}
    can serve as an argument to prove that 
    a given group $G \notin \mathcal{B}_f$ 
    for some $f \prec \mathfrak{i}$.    
 \end{rem} 
 
   Let $G$ be a group 
   $G= \langle A | R \rangle$ defined by 
   a finite set of generators $A$ and a 
   finite set of relators $R$. Let
   $S=A \cup A^{-1}$.  
   The Dehn function $D(n)$ of $G$ given 
   by $A$ and $R$ is defined as   
   $D(n) = \max \{\mathrm{Area}(w) |  
   w \in S^{\leqslant n} \wedge \pi(w)=e \}$, 
   where $\mathrm{Area}(w)$ is the minimal integer 
   $k$ for which 
   $w = \prod_{i=1}^{k} v_i r_i^{\pm 1} v_i^{-1}$,
   $r_i \in R$,
   in the free group $F(A)$.    
   Let us assume that 
   $G \in \mathcal{B}_f$ for some
   nonzero function $f \in  \mathfrak{F}$. 
   Theorem \ref{isoperimetric_theorem1} 
   and Corollary \ref{isoperimetric_corollary1}
   below extend the results we obtained in 
   \cite[Theorems~11 and 15]{measuring_closeness1}.    
   \begin{thm}        
   \label{isoperimetric_theorem1}      
      Assume that we are given two functions 
      $p,q \in \mathfrak{F}$ for which
      $p(n) \preceq D(n) \preceq q(n)$. Then 
      $p(n) \leqslant C n^2 q (K f(M n))$ for 
      all  $n \geqslant N$ for some 
      constants $C,K,M$ and $N$. 
      In particular, if 
      $p = q = n^d$ for some $d>2$, 
      then $n^{\frac{d-2}{d}} \preceq f$. 
      If $p = q = \mathfrak{e}$, 
      then $\mathfrak{i} \preceq f$. 
   \end{thm}
   \beginpf      
      Let $\psi: L \rightarrow G$ be 
      a Cayley automatic 
      representation of $G$ 
      such that for the function 
      $h(n) = \max \{d_A (\psi(w),\pi(w))|
       w \in L^{\leqslant n}\}$, $h \preceq f$.    
      Let $w = a_1 \dots a_n
      \in S^*$ be a word representing 
      the identity in $G$, where 
      $a_i \in S$. 
      For a given 
      $j=1,\dots,n-1$, we put 
      $g_j = a_1 \dots a_j$  
      and $g_0 = g_n = e$.  
      We first divide a loop given by 
      the word $w$ into  
      $n$ subloops as follows.
      For any $i=0,\dots,n-1$ 
      let 
      $u_i \in S^*$ 
      be the following concatenation of words:
      $u_i = \eta_i \xi_i a_{i+1} 
      \xi_{i+1}^{R} \eta_{i+1}^{R}$, 
      where $\eta_i = \psi^{-1}(g_i)$, 
      $\xi_i$ is some fixed word traversing 
      a shortest path from $\pi(\eta_i)$ 
      to $g_i$, 
      $\xi_{i+1}^{R}$ and $\eta_{i+1}^{R}$ 
      are the inverses of $\xi_{i+1}$ and 
      $\eta_{i+1}$, respectively; 
      e.g., if $\xi = abbc^{-1}a^{-1}$, then 
      $\xi^R = acb^{-1}b^{-1}a^{-1}$. 
      Clearly $\pi(u_i)=e$,
      so we obtain a loop. 
      
      By the bounded difference lemma 
      (see, e.g.,  
      \cite[Lemma~14.1]{KKM11}),      
      the length of each string $\eta_i$
      is bounded by $Cn$ for some constant $C$. 
      Then each of the subloops 
      given by $u_i$, $i=0,\dots,n-1$ 
      we divide into 
      at most $Cn$ smaller subloops
      as follows. 
      For every $1 \leqslant 
      j \leqslant \max\{|\eta_i|,|\eta_{i+1}|\}$  
      we construct a loop starting 
      at the point 
      $\widehat{\eta_i}(j-1)$ as follows. 
      For $1 \leqslant 
      j < \max\{|\eta_i|,|\eta_{i+1}|\}$
      the loop defined
      by the word 
      $v_{ij}= p_{ij}\zeta_{ij} 
       p_{(i+1)j}^R \zeta_{i(j-1)}^{R}$, where
       $p_{ij}$ is the string
       for which $\eta_i(j) = \eta_i (j-1) p_{ij}$ 
       (so $p_{ij}$ is either a single--letter
        string or the empty string)       
       and
       $\zeta_{ij}$ is some word traversing 
       a shortest path from $\widehat{\eta_i}(j)$ 
       to $\widehat{\eta_{i+1}}(j)$; 
       clearly, the length of this loop 
       is bounded by 
       $s(j)+s(j-1)+2$, 
       where $s$ is the function given 
       \eqref{def_function_f_1}.  
       For $j = \max \{|u_i|,|u_{i+1}|\}$ 
       the loop is defined by the  
       word $v_{ij} = p_{ij}\xi_i a_{i+1} 
       \xi_{i+1}^R p_{(i+1)j}^R 
       \zeta_{i(j-1)}^{R}$; 
       the length of this loop 
       is bounded by 
       $(2 h(j)+1)+s(j-1)+2$. 
       Let $\ell'(k)= \max\{2s(k)+2,
       2h(k)+s(k)+3\}$ and 
       $\ell (k) = \ell' (Ck)$.       
       So, the length of each of these
       smaller subloops is bounded 
       by $\ell(n)= \ell'(Cn)$. 
       By the inequalities
       $h \preceq f$ and 
       $s \preceq f$       
       (see Theorem \ref{fellow_theorem1}), 
       we have $\ell \preceq f$. 
       The total number of these smaller subloops
       is at most $Cn^2$. Thus we obtain
       the inequality 
       $D(n) \leqslant Cn^2 D(\ell(n))$. Therefore,
       $D(n) \preceq n^2 D(\ell(n))$.

      From the inequalities
      $D(n) \preceq n^2 D(\ell(n))$, 
      $\ell \preceq f$      
      and      
      $p(n)\preceq D(n) \preceq q(n)$        
      we obtain that: 
      $p(n) \leqslant C_1 D(C_2 n) 
      \leqslant C_3 n^2 D(\ell (C_4 n)) 
      \leqslant C n^2 q(C_5 \ell (C_4 n)))
      \leqslant C n^2 q(K f(M n))$ for
      all $n \geqslant N$ for some constants 
      $C,K,M,N$ and $C_i,i=1,\dots,5$.   
      If $p = q = n^d$, then 
      $n^d \leqslant C n^2 (K f(Mn))^d$ 
      for all $n \geqslant N$. Therefore, 
      $n^\frac{d-2}{d} \leqslant 
       C^\frac{1}{d} K f(Mn)$ for all 
      $n \geqslant N$, i.e., 
      $n^\frac{d-2}{d} \preceq f$. 
      If $p=q=\mathfrak{e}$, then 
      $\exp(n) \leqslant C n^2 \exp(Kf(Mn))$ 
      for all $n \geqslant N$. 
      Therefore, $n \leqslant \log C + 
      2 \log n + K f(Mn)$ for all 
      $n \geqslant N$, which implies that
      $\mathfrak{i} \preceq f$.             
   \endpf

  \begin{cor} 
  \label{isoperimetric_corollary1}     
     For a given function $f \in \mathfrak{F}$ 
     we have: 
     \begin{itemize}
       \item{if the Baumslag--Solitar group 
       $BS(p,q) \in \mathcal{B}_f$ for some 
       $q > p \geqslant 1$, then 
       $\mathfrak{i} \preceq f$;}
       \item{if the Heisenberg group 
       $\mathcal{H}_3(\mathbb{Z}) \in 
       \mathcal{B}_f$, then $\sqrt[3]{n} 
       \preceq f$;}
       \item{if the group 
       $\mathbb{Z}^2 \rtimes_A \mathbb{Z} 
       \in \mathcal{B}_f$ for a matrix 
       $A \in GL(2,\mathbb{Z})$ with 
       two real eigenvalues not equal to 
       $\pm 1$, then $\mathfrak{i} \preceq f$.} 
     \end{itemize}
  \end{cor}   
  \beginpf
     This follows from 
     Theorem \ref{isoperimetric_theorem1} 
     and the facts that 
     for the groups 
     $BS(p,q)$, $1\leqslant p<q$, 
     $\mathcal{H}_3(\mathbb{Z})$ and 
     $\mathbb{Z}^2 \rtimes_A \mathbb{Z}$,
     for a matrix 
     $A \in GL(2,\mathbb{Z})$ with 
     two real eigenvalues not equal to $\pm 1$,     
     the Dehn functions are exponential, 
     cubic and exponential, respectively 
     (see \cite{BrisonPittet94} and, e.g., 
     \cite[\S7.4--\S8.1]{Epsteinbook}).      
   \endpf
    \begin{rem}    
      We recall that the groups 
      $\mathbb{Z}^2 \rtimes_A \mathbb{Z}$ 
      are the fundamental groups of 
      $3$--manifolds which are 
      $2$--dimensional torus bundles 
      over the circle. The Heisenberg group 
      $\mathcal{H}_3(\mathbb{Z})$ is isomorphic
      $\mathbb{Z}^2 \rtimes_A \mathbb{Z}$ 
      for some unipotent matrix $A$; see also 
      Section \ref{nonautomaticitysection}.       
    \end{rem}    
    \begin{rem}
       The examples of Dehn functions for
       Cayley automatic groups, which 
       are known to us,
       are quadratic (e.g, for 
       the higher Heisenberg
       groups $\mathcal{H}_{2k+1}(\mathbb{Z}),
       k>1$),
        cubic 
       (e.g., for the Heisenberg group
       $\mathcal{H}_3(\mathbb{Z})$), 
       $n^d$ for any integer $d>3$ 
       (e.g., for some semidirect products 
        $\mathbb{Z}^m \rtimes_A \mathbb{Z}$,
        see 
        \cite{BridsonGersten96,BrisonPittet94}), 
        and 
        the exponential function $\mathfrak{e}$ 
        (e.g., for the Baumslag--Solitar groups
        $BS(p,q)$, $1 \leqslant p <q$). 
    \end{rem}

 \section{{\bf Finite Extensions, 
 Direct Products,\\ Free Products}} 
 \label{extfreeproddirprod}

 In this section we show that classes 
 $\mathcal{B}_f$ are invariant  with 
 respect to taking finite extension, 
 direct product and free product. For the
 latter case we require that 
 $f$ satisfies the inequality 
 $f(x)+f(y) \leqslant f(x+y)$ for 
 all $x,y \geqslant n_0$, where 
 $n_0$ is some constant.  
 Let $H$ be a subgroup of finite 
 index in a f.g. group $G$. 
 It is known that if $H$ is automatic, then
 $G$ is automatic.
 Moreover, by \cite[Theorem~10.1]{KKM11}, if 
 $H$ is Cayley automatic, then $G$ is Cayley
 automatic\footnote{A complete analog of 
 \cite[Theorem~4.1.4]{Epsteinbook} 
 for automatic groups, claiming that 
 a subgroup $H$ of finite index of a group $G$ is 
 automatic iff $G$ is automatic, is not known 
 for Cayley automatic groups. We remark that 
 in the original \cite[Theorem~10.1]{KKM11}
 the assumption that $H$ is a normal subgroup 
 of $G$ can be omitted; 
 see, e.g., 
 \cite[Theorem~2.2.4]{BerdinskyCSThesis}.}.       
 \begin{thm}
 \label{finiteindex_theorem1}    
    Let $H$ be a subgroup of finite index of a 
    group $G$. If $H \in \mathcal{B}_f$, 
    then $G \in \mathcal{B}_f$.  
 \end{thm}
 \beginpf
    Let us fix a finite set of generators
    of $H$: $A_1 = \{h_1,\dots,h_n\}$, and 
    a set of unique representatives of the right 
    cosets $Hg$ of the subgroup $H$ in 
    $G$, where $g \notin H$: $A_2 = 
    \{k_1,\dots,k_m\}$. 
    We put $S_1 = A_1 \cup A_1^{-1}$. 
    Since $H \in \mathcal{B}_f$, there exist 
    a Cayley automatic representation 
    $\psi_1: L_1 \rightarrow H$,     
    $L_1 \subseteq S_1^*$
    such that, for the function
    $h_1(n)= \max \{ d_{A_1} (\pi(u),\psi_1(u)) 
    \,|\, u \in L_1 ^{\leqslant n}\}$,
    $h_1(n) \preceq f(n)$. 
    Let $L_2$ be the finite language consisting 
    of $m$ single--letter strings
    $k_1,\dots,k_m$ and the empty
    sting $\epsilon$. We put $\psi_2$ to be the 
    natural embedding of these strings into 
    the group $G$: a string $k_i$ maps 
    to the group 
    element $k_i$ and the empty 
    string $\epsilon$
    maps to the identity of the group $G$.   
    We put $L$ to be the concatenation
    of $L_1$ and $L_2$. Clearly, 
    $L \subseteq S^*$, where 
    $S = A \cup A^{-1}$ and 
    $A = A_1 \cup A_2$.   
    Now, we define the map 
    $\psi: L \rightarrow G$ as follows. 
    Let $w = uv \in L$, where $u \in L_1$ 
    and $v \in L_2$. We put 
    $\psi(w):=\psi_1(u)\psi_2 (v)$. 
    It is easy to verify that the constructed 
    map $\psi$ is a Cayley 
    automatic representation 
    of the group $G$ 
    (see \cite[Theorem~10.1]{KKM11}). 
    Furthermore, 
    $d_A (\pi(w),\psi(w)) \leqslant 
    d_A (\pi(w),\pi(u)) + d_A (\pi(u),\psi(u))
    + d_A (\psi(u),\psi(w)) \leqslant 
    1 + h_1(|u|) + 1 \leqslant h_1 (|w|)+2$. 
    This immediately implies that      
    for the function 
    $h(n)= \max 
    \{d_A(\pi(w),\psi(w))\,|\,
    w \in L^{\leqslant n}\}$, 
    $h \preceq h_1$. Therefore, 
    $h \preceq f$.       
 \endpf
  
  It is known that the direct product of two 
 automatic groups is automatic. 
 The direct product of Cayley automatic 
 groups is also Cayley automatic 
 \cite[Corollary~10.4]{KKM11}. 
 \begin{thm}
 \label{directprodthm1}    
    If $G_1,G_2 \in \mathcal{B}_f$, then 
    $G_1 \times G_2\in \mathcal{B}_f$.  
 \end{thm}
 \beginpf
    Let $A_1$ and $A_2$ be some sets 
    of generators 
    of the groups $G_1$ and $G_2$ for which
    $A_1 \cap A_2 = \varnothing$; 
    we put $S_1 = A_1 \cup A_1^{-1}$ and 
    $S_2 = A_2 \cup A_2^{-1}$.      
    Since $G_1,G_2 \in \mathcal{B}_f$, there
    exist Cayley automatic representation 
    $\psi_1: L_1 \rightarrow G_1$ and 
    $\psi_2: L_2 \rightarrow G_2$ for which 
    the functions $h_1(n)= \max \{ 
    d_{A_1}(\pi(w),\psi_1(w))
    \,|\,w \in L_1 ^{\leqslant n} \}$ and 
    $h_2 (n)= \max \{ 
    d_{A_2}(\pi(w),\psi_2(w))\,|\,w \in 
    L_2 ^{\leqslant n}\}$ satisfy the 
    inequalities $h_1 \preceq f$ and
    $h_2 \preceq f$, where 
    $L_1 \subseteq S_1 ^*$ and 
    $L_2 \subseteq S_2 ^*$.     
    
    Let $L = L_1 L_2$. We construct 
    the map $\psi: L \rightarrow 
    G_1 \times G_2$ as follows.     
    For a given $w= uv$, where $u \in L_1$ 
    and $v \in L_2$, we put 
    $\psi(w) = (\psi_1 (u),\psi_2 (v))
     \in G_1 \times G_2$. 
     It is easy to verify that the constructed  
    map $\psi$ provides a Cayley automatic 
    representation of $G_1 \times G_2$.  
    The groups $G_1$ and 
    $G_2$ are naturally embedded in  
    $G_1 \times G_2$, so we have
    $\pi(w)=\pi(u) \pi(v) = 
     (\pi(u),\pi(v)) \in G_1 \times  G_2$.        
    Therefore, 
    $d_A (\pi(w),\psi(w))  
          \leqslant 
          d_A (\pi(u),\psi_1(u))+
          d_A (\pi(v),\psi_2(v))\leqslant
          h_1 (|u|) + h_2 (|v|) \leqslant
          h_1(|w|) + h_2(|w|)=s(|w|)$, 
          where $s(n)= h_1 (n)+h_2(n)$ for
          all $n \in \mathrm{dom}\,h_1 \cap 
                 \mathrm{dom}\,h_2$.  
     Clearly, the inequalities $h_1 \preceq f$ 
              and $h_2 \preceq f$ 
              imply that $s \preceq f$. 
     Therefore, for the function 
     $h(n)=\max\{d_A(\pi(w),\psi(w))|
           w \in L^{\leqslant n}\}$, we have 
     $h \preceq f$. 
 \endpf 
 
  It is known that the free product 
  of automatic groups is automatic. 
  Therefore, if $G_1,G_2 \in \mathcal{B}_d$, 
  then $G_1 \star G_2 \in \mathcal{B}_d$, where 
  $d$ is a bounded function (recall that
  in this case,    
  by \cite[Theorem~8]{measuring_closeness1},  
  $\mathcal{B}_d$ is 
  the class of automatic groups). 
  Moreover, the free product of Cayley 
  automatic groups is Cayley automatic 
  \cite[Theorem~10.8]{KKM11}.        
  In the following
  theorem we consider the case when 
  $G_1,G_2 \in \mathcal{B}_f$ for some
  unbounded function $f \in \mathfrak{F}$. 
 \begin{thm} 
 \label{freeprodtheorem1}    
    Let $f \in \mathfrak{F}$ be a function
    for which $f(x)+ f(y) \leqslant f(x+y)$
    for all $x,y \geqslant n_0$, where 
    $n_0$ is a constant.     
    If $G_1,G_2 \in \mathcal{B}_f$, 
    then $G_1 \star G_2 \in \mathcal{B}_f$. 
 \end{thm}
 \beginpf 
    For initial settings 
    we use the same notation as  
    in the first 
    paragraph  of the proof of 
    Theorem \ref{directprodthm1}.     
    Without loss of generality we may assume 
    that the empty word $\epsilon \in L_1,L_2$,  
    and $\psi_1 (\epsilon)$ and 
    $\psi_2 (\epsilon)$ are the identities 
    in the groups $G_1$ and $G_2$,
     respectively. 
    We put
    $L_1' = L_1 \setminus \{ \epsilon \}$ 
    and 
    $L_2'= L_2 \setminus \{\epsilon \}$.     
    Let $A = A_1 \cup A_2$.   
    Let $L$ be defined by the following regular 
    expression 
    $L = (L_1' L_2')^* \vee
         (L_1' L_2')^* L_1' \vee 
         (L_2' L_1')^* \vee         
         (L_2' L_1')^* L_2' 
         \vee \epsilon$. 
    That is, $L$ is the regular language 
    consisting of the empty string 
    $\epsilon$ and the strings of the form 
    $u_1 \dots u_k$, where each 
    substring $u_i$, $i=1,\dots,k$ either 
    $u_i \in L_1 '$ or 
    $u_i \in L_2 '$,     
    and no consecutive strings 
    $u_i,u_{i+1}$ are 
    elements of the same language
     $L_1'$ or $L_2'$. 
    Let us construct the map 
    $\psi: L \rightarrow G_1 \star G_2$ 
    as follows: 
    $\psi(\epsilon)=e$ and 
    $\psi (u_1 \dots u_k)=\psi (u_1) 
    \dots \psi(u_k)$, where 
    for each $u_i$, $i=1,\dots,k$, 
    $\psi(u_i)=\psi_1 (u_i)$ or
    $\psi(u_i)=\psi_2 (u_i)$ if 
    $u_i \in L_1'$ or $u_i \in L_2'$,
     respectively. 
    It is easy to verify that the constructed  
    map $\psi$ provides a Cayley automatic 
    representation of $G_1 \star G_2$      
    (see also \cite[Theorem~10.8]{KKM11}). 
    
    Now, let $w = u_1 \dots u_k \in L$. Then, 
    $d_A (\pi(w),\psi(w)) \leqslant 
     d_A (\pi(w)) + d_A (\psi(w)) \leqslant 
     |w| + \sum_{i=1}^k d_A (\psi(u_i))$. 
   For each $u_i$, $i=1,\dots,k$, 
   we have $d_A(\psi(u_i)) \leqslant 
   d_A (\pi(u_i))+d_A(\pi(u_i),\psi(u_i)) 
   \leqslant |u_i| + K f (M|u_i|)$, 
   if $|u_i| \geqslant N$ for some positive  
   integer constants
   $K,M$ and $N$; here we also assume that 
   $MN \geqslant n_0$. For all $|u_i| < N$ we 
   can bound $d_A (\psi(u_i))$ from above 
   by some constant $C$ 
   since there exist only finitely 
   many such $u_i$; we also assume that 
   $C \geqslant 1$. 
   Therefore, by the assumption that 
   $f(x)+f(y) \leqslant f(x + y)$ for all 
   $x,y \geqslant n_0$,   
   we obtain  
   $\sum_{i=1}^k d_A (\psi(u_i)) \leqslant 
    C |w| + K  f(M|w|)$.    
   Thus, $d_A (\pi(w),\psi(w)) \leqslant
   (C+1) |w| + K f(M|w|)$ for all  
   $w \in L$. 
   We note that the inequality 
   $f(x)+f(y)\leqslant f(x+y)$ 
   for all $x,y \geqslant n_0$ implies 
   that $\mathfrak{i} \preceq f$, unless
   $f$ is identically equal to zero.   
   So, for the function 
   $h (n) = \max \{d_A (\pi(w),\psi(w)) | 
    w \in L^{\leqslant n} \}$, we 
    have $h \preceq f$.      
 \endpf
 \begin{cor}
  If
    $G_1,G_2 \in \mathcal{B}_{\mathfrak{i}}$ 
    or $G_1,G_2 \in \mathcal{B}_{\mathfrak{e}}$, 
    then $G_1 \star G_2$ is also 
    in the class $\mathcal{B}_{\mathfrak{i}}$
    or $\mathcal{B}_{\mathfrak{e}}$, respectively.
 \end{cor}
 \beginpf
    It is enough to notice that for the functions 
    $f = \mathfrak{i}$ and 
    $f = \mathfrak{e}$,    
    the inequality 
    $f(x)+f(y) \leqslant f(x+y)$ holds 
    for all $x,y \geqslant 1$.
 \endpf

 \section{{\bf Nilpotent Groups and
    Fundamental\\ Groups of $n$--dimensional 
    Torus Bundles over\\ The Circle}}
 \label{nilpotent_section} 
  
   In this section we show that some 
   classes of nilpotent groups    
   and the fundamental groups of 
   $n$--dimensional torus bundles 
   over the circle are in the class 
   $\mathcal{B}_\mathfrak{e}$.
   In the second half of the section we  
   address the problem of finding 
   sharp lower bounds of 
   the function \eqref{def_func_h_1} 
   for virtually nilpotent groups.  
   Before we proceed with the main result 
   of the section let us prove the following 
   technical lemma which is needed,
   in particular, for the proof 
   of Theorem \ref{nilpotentsec_thm1}.    
     Let $\varphi:L_\Sigma 
     \rightarrow G$ be a Cayley 
     automatic representation of $G$, 
     where $L_\Sigma \subseteq \Sigma^*$ now is a 
     regular language over some alphabet 
     $\Sigma^*$ (here we do not assume
     that $\Sigma = S$). We denote by 
     $h_\varphi$ the function
     $h_\varphi (n) = \max \{ 
     d_A (\varphi(w)) | w \in L ^{\leqslant n}\}$.

  \begin{lem}
  \label{facts_lemma1}  
     Suppose that $h_\varphi \preceq f$ 
     for some function $f \in  \mathfrak{F}$. 
     Then $G \in \mathcal{B}_{\tilde{f}}$, 
     where $\tilde{f} = f+\mathfrak{i}$.    
  \end{lem} 
  \beginpf
     For every $\sigma \in \Sigma$ let us
     choose a string $w_\sigma \in S^*$
     such that the lengths $|w_\sigma|$ 
     are equal to some constant $\ell$ 
     for all $\sigma \in \Sigma$. 
     Then we define a 
     monoid homomorphism 
     $\xi : \Sigma^* \rightarrow S^*$
     as follows: 
     $\xi (\sigma_1 \dots \sigma_k ) = 
      w_{\sigma_1} \dots w_{\sigma_k}$. 
     We define $L = \xi (L_\Sigma)$ and  
     $\psi = \varphi \circ \xi^{-1}: 
      L \rightarrow G$. Clearly, 
      $\psi: L \rightarrow G$ is a Cayley 
      automatic representation of $G$. 
     Moreover, for any $w \in L$ we have 
     $d_A(\pi(w),\psi(w)) \leqslant 
      d_A(\pi(w))+d_A(\psi(w)) = |w| + 
      d_A(\varphi \circ \xi^{-1}(w)) 
      \leqslant |w| + h_\varphi (|\xi^{-1}(w)|) 
      = |w| + h_\varphi 
      \left(\frac{1}{\ell} |w| \right) 
      \leqslant |w| + h_\varphi (|w|)$. 
      Therefore, for the function 
      $h(n)= \max \{d_A (\pi(w),\psi(w)) | 
             w \in L ^{\leqslant n}\}$, we
      clearly have $h \preceq  
            \tilde{f}$.       
      \endpf   
 
  \begin{thm}
  \label{nilpotentsec_thm1}     
     The following groups are all in the class 
     $\mathcal{B}_\mathfrak{e}$: 
     \begin{itemize}       
        \item{fundamental groups 
        of $n$--dimensional torus bundles 
        over the circle     
        
        $\mathbb{Z}^n \rtimes _A 
        \mathbb{Z}$,} 
        \item{unitriangular matrices
        $UT_n(\mathbb{Z})$,}    
        \item{f.g. nilpotent groups 
        of nilpotency class $2$.}     
     \end{itemize}
  \end{thm}  
  \beginpf
     Let $\beta$ be a representation 
     of $\mathbb{Z}$ for which every 
     $z \in \mathbb{Z}$ is represented 
     as a signed binary number.  
     Let $\gamma$ be a representation of
     $\mathbb{Z}$  for which 
     every $y \in \mathbb{Z}$ is represented 
     as the concatenation of $|y|$ 
     identical single--letter strings; 
     for positive and negative integers we 
     use different letters.      
     See also the representation of the 
     Heisenberg group 
     $\mathcal{H}_3 (\mathbb{Z})$ 
     that we constructed
     in \cite[Section~6]{measuring_closeness1}.          
     For any given 
     $\overline{z} = 
     (z_1,\dots,z_n) \in \mathbb{Z}^n$ 
     we represent it as the convolution 
     $v = 
     w_1 \otimes \dots \otimes w_n$, where
     $w_i = \beta^{-1} (z_i)$, 
     $i=1, \dots, n$.     
     Then we represent      
     an element  
     $g= (y,\overline{z}) \in 
     \mathbb{Z}^n \rtimes_A \mathbb{Z}$
     as the concatenation 
     $w= u v$, where $u = \gamma^{-1}(y)$. 
     By \cite[Theorem~10.3]{KKM11}, 
     it provides a Cayley automatic 
     representation $\varphi$ of 
     $\mathbb{Z}^n \rtimes_A \mathbb{Z}$. 
     In the group  
     $\mathbb{Z}^n \rtimes_A \mathbb{Z}$
     the element $g = (y,\overline{z})$ is  
     equal to the product 
     $g = (y,\overline{0}) \cdot (0,\overline{z})$, 
     where $0$ is the identity of $\mathbb{Z}$ and
     $\overline{0}$ is the identity of 
     $\mathbb{Z}^n$. It is easy to see now 
     that the condition of Lemma 
     \ref{facts_lemma1} is satisfied 
     for the representation 
     $\varphi$,      
     the function $f = \mathfrak{e}$
     and      
     a natural set of generators 
     $(1,\overline{0})$ and 
     $(0,\overline{e}_i)$, 
     $i=1,\dots,n$, where
     $\overline{e}_i \in \mathbb{Z}^n$ 
     has the $j$th element equal to 
     $\delta_{ij}$, $j=1,\dots,n$.       
     Therefore, $\mathbb{Z}^n \rtimes _A 
     \mathbb{Z} \in \mathcal{B}_\mathfrak{e}$. 
     
     Any element $g$ of 
     the unitriangular matrix group     
     $UT_n (\mathbb{Z})$ 
     is given by a $n \times n$ matrix $M$ 
     with all elements below the main 
     diagonal equal to $0$ and 
     all elements of the main diagonal equal
     to $1$.       
     Let $m_{ij} \in \mathbb{Z}$, $i<j$ 
     be the element of $M$ 
     in row $i$ and column $j$.  
     We denote by $t_{ij} \in UT_n (\mathbb{Z})$   
     the transvection given by a $n \times n$
     matrix with all elements on the 
     main diagonal and the element in 
     row $i$ and column $j$ equal to $1$ and
     all other elements equal to $0$.     
     In the group $UT_n(\mathbb{Z})$ 
     the element $g$ is equal to the product
     of transvections 
     $g = t_{1n}^{m_{1n}} \dots  
      t_{(n-1)n}^{m_{(n-1)n}} 
      \dots      
      t_{13}^{m_{13}} t_{23}^{m_{23}} 
      t_{12}^{m_{12}}$.
     We represent $g$ as the convolution 
     $s_{12} \otimes \dots 
             \otimes s_{(n-1)n}$, 
     where $s_{ij}= \beta^{-1}(m_{ij})$, 
     $1 \leqslant i < j \leqslant n$. 
     Clearly, the condition of Lemma 
     \ref{facts_lemma1} is satisfied 
     for this representation, the function 
     $f= \mathfrak{e}$ and the set of 
     generators 
     $\{ t_{ij} | 
      1 \leqslant i < j \leqslant n \}$. 
      Therefore, $UT_n (\mathbb{Z}) \in
      \mathcal{B}_\mathfrak{e}$.
      
     It is known that 
     for every f.g. nilpotent 
     group its torsion subgroup is finite. 
     Moreover, every f.g. nilpotent group 
     is residually finite. Therefore, 
     every f.g. nilpotent group has a 
     torsion--free subgroup of finite index.
     So, by Theorem \ref{finiteindex_theorem1}, 
     it is enough for us to show that any  
     given torsion--free f.g. nilpotent group $G$ 
     of nilpotency class $2$ is in  
     $\mathcal{B}_\mathfrak{e}$.      
     In \cite[Theorem~12.4]{KKM11} 
     the authors used Mal'cev coordinates 
     to construct Cayley automatic representation
     of the group $G$. Below we use their 
     representation to show that 
     $G \in \mathcal{B}_\mathfrak{e}$.    
     Let $ \overline{a}= 
     (a_1,\dots,a_n) \in G^n$ be any Mal'cev 
     basis for $G$ associated with the upper 
     central series of $G$. We recall that 
     the factors of the upper central series 
     of a torsion--free nilpotent group are 
     torsion--free.   
     So, for any given $g \in G$, 
     we have a unique presentation
     of $g$ in $G$ as a product: 
     $g = a_1^{k_1} \dots 
               a_n^{k_n}$, 
     where $(k_1,\dots,k_n) \in 
     \mathbb{Z}^n$ is a tuple of
      the Mal'cev coordinates of $g$ 
      with respect to the basis 
      $\overline{a}$. 
     We represent $g$ as the convolution 
     $s_1 \otimes \dots  \otimes s_n$, where 
     $s_i = \beta^{-1}(k_i)$, $i=1,\dots,n$. 
     The condition of Lemma \ref{facts_lemma1} 
     is satisfied for this representation, 
     the function $f = \mathfrak{e}$ and 
     the set of generators 
     $\{a_1,\dots,a_n\}$. Thus, 
     $G \in \mathcal{B}_\mathfrak{e}$.  
  \endpf
   
    Can any of the groups from Theorem 
    \ref{nilpotentsec_thm1} be in the class 
    $\mathcal{B}_f$ for some 
    $f \prec \mathfrak{e}$? 
    The greatest lower bound for the function 
    $f$ that we can obtain from Theorem 
    \ref{isoperimetric_theorem1} is 
    $\mathfrak{i}$, see, e.g., Corollary 
    \ref{isoperimetric_corollary1}. 
    However, for some groups, e.g.
    the higher Heisenberg
    groups $\mathcal{H}_{2k+1}$, $k>1$, 
    Theorem \ref{isoperimetric_theorem1} 
    does not give any lower bound 
    (recall that they are nilpotent groups 
    of nilpotency class $2$ and their 
    Dehn functions
    are quadratic).         
    Thurston proved that automatic 
    nilpotent groups must be virtually abelian
    (see, e.g., \cite[Theorem~8.2.8]{Epsteinbook}).  
    So, by 
    \cite[Theorem~8]{measuring_closeness1}, 
    for any class $\mathcal{B}_f$ 
    containing  a Cayley automatic nilpotent group 
    (which is not virtually abelian) 
    the function $f$ must be unbounded.  
    Moreover, while for the Baumslag--Solitar 
    groups $BS(p,q)$, $q>p \geqslant 1$ and 
    the lamplighter group 
    $\mathbb{Z}_2 \wr \mathbb{Z}$ we obtain 
    the sharp lower bounds
    \cite[Theorem~11 and 13]{measuring_closeness1}, 
    we do not know whether the lower bounds,
    which we can obtain from Theorem 
    \ref{isoperimetric_theorem1}
    for other groups mentioned in this paper, 
    are sharp. To address this issue we make 
    a simple   
    observation in Theorem
    \ref{nilpotentsec_starredlangthm1} 
    that might, potentially, be useful 
    in the search for
    the sharp lower bounds 
    for virtually nilpotent groups.
    Furthermore, 
    in Theorem \ref{nonautth1} we show
    that, for the 
    Heisenberg group 
    $\mathcal{H}_3(\mathbb{Z})$,
    the exponential function 
    $\mathfrak{e}$ is    
    a lower bound of the function 
    \eqref{def_func_h_1},   
    if one puts some
    additional constraints on a Cayley 
    automatic representation $\psi$.    
   We recall that a regular language 
   $L$ is called 
   simply starred if      
   a regular expression for 
     $L$ is of the form: 
     $
        R_1 \vee \dots  \vee R_I,       
     $
     where $R_i = 
            v_{i,0} u_{i,1}^* v_{i,1} 
            \dots  v_{i, P_{i-1}} 
            u_{i,P_i}^* v_{i,P_i}$ 
            for $i=1,\dots,I$.   
   We have the following proposition.      
  \begin{prop}[polynomial growth condition] 
  \label{reglanpolgrowthcond}   
     A regular language $L$ has polynomial 
     growth if it is simply starred and 
     exponential growth otherwise.      
  \end{prop}  
  \beginpf
     For the proof see, e.g., 
     \cite[Theorem~8.2.8]{Epsteinbook}.  
  \endpf
  
  Let $\psi: L \rightarrow G$ be a 
  Cayley automatic representation 
     of a virtually nilpotent group 
     $G$; as usual,  
     $L \subseteq (A \cup A^{-1})^*$ for some 
     set of generators $A \subset G$. 
     Let $h$ be the function 
     defined by \eqref{def_func_h_1} 
     corresponding to 
     the representation $\psi$. 
  \begin{thm}      
  \label{nilpotentsec_starredlangthm1}     
     Suppose that   
     $h \preceq p$ for some 
     polynomial $p$.  
     Then the language $L$
     is simply starred.
  \end{thm}
  \beginpf
      For any given $w \in L^{\leqslant n}$ 
      we have
      $d_A (\psi(w)) \leqslant 
       d_A (\pi(w)) + d_A (\pi(w),\psi(w))$ 
       $\leqslant n + h (n)$. 
       Therefore, since $h \preceq p$, 
       there exists a  polynomial $q$ 
       for which $\psi(w)$ must be 
       in the ball $B_{q(n)}\subset G$ 
       of radius $q(n)$. 
       Recall that a growth function of any 
       virtually nilpotent group 
       is bounded by a polynomial. Therefore, 
       the cardinality of  
       $B_{q(n)}$ must be bounded by
       $r(n)$ for some 
       polynomial $r$ so the cardinality of 
       the set
       $L^{\leqslant n}$. By    
       Proposition \ref{reglanpolgrowthcond} 
       we obtain the statement of the theorem.     
  \endpf
  
 \section{\bf In The Search for 
 Alternative Approaches to 
 Proving Nonautomaticity } 
 \label{nonautomaticitysection}
      
    In this section we focus on the problem 
    of finding a sharp lower bound
    of the function \eqref{def_func_h_1} for 
    the Heisenberg group  
    $\mathcal{H}_3 (\mathbb{Z})$.
    Another motivation of 
    this section is to propose 
    alternative methods 
    for proving nonautomaticity of groups.     
    Clearly, if a group 
    $G \notin \mathcal{B}_f$ 
    for some  
    function 
    $f \in \mathfrak{F}$, 
    then $G$ is not automatic. 
    We already know two ways 
    to show that a group is not in 
    a class $\mathcal{B}_f$ if 
    $f \prec f_0$ for some nonzero 
    function $f_0$ (see 
    Theorem \ref{isoperimetric_theorem1} and 
    the proof that the lamplighter group 
    is not in the class $\mathcal{B}_f$ for
    any $f \prec \mathfrak{i}$ 
    \cite[Theorem~13]{measuring_closeness1}). 
    In the first 
    approach we 
    use the Dehn function 
    (when it grows faster than the quadratic 
     function), while in the second 
     approach we implicitly 
     use a fact that the 
     lamplighter group 
     is not finitely presented. 
     However, in both cases 
     one straightforwardly gets nonautomaticity 
     by \cite[Theorem~2.3.12]{Epsteinbook}.
     Is there any alternative method to show
     that a given group 
     $G$ is not in $\mathcal{B}_f$ for 
     some $f \in \mathfrak{F}$? 
     Such a method could potentially  provide  
     a new way to prove nonautomaticity. 
     In this part we make a first tiny step 
     in this direction focusing on the 
     Heisenberg group 
     $\mathcal{H}_3 (\mathbb{Z})$. 
    
     It was first noticed by S\'{e}nizergues
     that the Heisenberg group is not
     automatic, 
     but its Cayley graph is FA--presentable;
     also, 
     it was one of the first examples of 
     such groups.            
     Another motivation to focus on 
     $\mathcal{H}_3(\mathbb{Z})$ is the 
     "Heisenberg alternative"  
     -- each f.g. group $G$ of 
     polynomial growth is either virtually 
     abelian or 
     $\mathcal{H}_3 (\mathbb{Z})$ 
     can be embedded into $G$. 
     In \cite{NiesThomas08},
     Nies and Thomas used this alternative  
     to give a new proof of the 
     theorem that every f.g. 
     FA--presentable group is virtually abelian; 
     this was first proved by 
     Oliver and Thomas in \cite{OliverThomas05}.          
     We recall that 
     $\mathcal{H}_3 (\mathbb{Z})$
     is the group of all matrices of the form: 
     $ \left(
      \begin{array}{ccc}
       1 & x & z \\
       0 & 1 & y \\
       0 & 0 & 1
     \end{array}
     \right),$ 
     where $x,y$ and $z$ are integers; so,  
     every element 
     $g \in \mathcal{H}_3 (\mathbb{Z})$
     corresponds to a triple $(x,y,z)$.      
     We denote by  
     $s,p$ and $q$ the group elements 
     corresponding to the triples 
     $(1,0,0)$, $(0,1,0)$, and 
     $(0,0,1)$, respectively. 
   If $g$ corresponds to a triple $(x,y,z)$, 
   then $gs,gp$ and $gq$ 
   correspond to the triples 
   $(x+1,y,z)$, $(x,y+1,x+z)$ and $(x,y,z+1)$, 
   respectively. We put 
   $A = \{e,s,p,q\}$ and 
   $S = \{e,s,p,q,s^{-1},p^{-1},
    q^{-1}\}$.

     It is straightforward to verify  
     that $\mathcal{H}_3 (\mathbb{Z})$ 
     is isomorphic 
     to the semidirect product 
     $\mathbb{Z}^2 \rtimes_T \mathbb{Z}$, 
     where  
     $T = \left(\begin{array}{cc}
      1 & 0 \\     
      1 & 1 
      \end{array}
      \right)$: an isomorphism is given 
      by the following mapping
      $(x,y,z) \mapsto 
       \left(y,\left[\begin{array}{c}
       x \\     
       z 
       \end{array}
      \right]\right)$.  
      We denote by $H$ the normal subgroup 
      of $\mathcal{H}_3 (\mathbb{Z})$
      generated by  $s$ and $q$, 
      and by $N$ the cyclic subgroup 
      of $\mathcal{H}_3 (\mathbb{Z})$ 
      generated by $p$. Clearly, 
      $H \cong \mathbb{Z}^2$, 
      $N \cong \mathbb{Z}$ and 
      $\mathcal{H}_3 (\mathbb{Z}) = NH$.    
      We denote by $\varphi$, $p_1$ and $p_2$
      the endomorphisms of the group $H$ given by 
      the matrices $T = \left(\begin{array}{cc}
      1 & 0 \\     
      1 & 1 
      \end{array}
      \right)$, $P_1= \left(\begin{array}{cc}
      1 & 0 \\     
      0 & 0 
      \end{array}
      \right)$ and 
      $P_2 = \left(\begin{array}{cc}
      0 & 0 \\     
      0 & 1 
      \end{array}
      \right),$ respectively. The endomorphisms 
      $p_1$ and $p_2$ are the projectors
      of $H$ on the cyclic subgroups  
      generated by $s$ and $q$, respectively. 
      We denote these subgroups by $H_1$ and 
      $H_2$:  
      $H_1 = p_1(H)= \langle s \rangle$ and 
      $H_2 = p_2(H)= \langle q \rangle$.    
      Let $\psi: L \rightarrow 
             \mathcal{H}_3 (\mathbb{Z})$ be 
             a Cayley automatic 
             representation of 
             $\mathcal{H}_3(\mathbb{Z})$, 
        where $L \subseteq S^*$. 
      We denote by $L_H$ the language
                   $L_H = \psi^{-1}(H)
                   \subset L$ and 
                   by $w_0$ the string 
                   $w_0 = \psi^{-1}(e)$,
                   where $e$ is the identity 
                   of the group 
                   $\mathcal{H}_3(\mathbb{Z})$.      
      Let 
      $R_\varphi = 
      \{\langle w, 
      \psi^{-1} \circ \varphi \circ \psi (w)
      \rangle \,|\, w \in L_H \}, 
      R_{p_1}= \{ \langle w, 
      \psi^{-1} \circ p_1 \circ \psi (w)
      \rangle \, | \, w \in L_H \},
      R_{p_2} = \{ \langle w, 
      \psi^{-1} \circ p_2 \circ \psi (w) 
      \rangle \,| \, w \in L_H \}       
      \subset L_H \times L_H$ 
      be the binary relations on $L_H$  
      defined by the 
      endomorphisms $\varphi,p_1,p_2$, 
      respectively, and the 
      Cayley automatic representation $\psi$. 
      For a given binary relation 
      $R \subseteq S^* \times S^*$,     
      we denote by $L_H \triangleleft R$ 
      and $R \triangleright L_H$
      the left-- 
      and right--restrictions of 
      $R$ on $L_H$: $L_H \triangleleft R= 
      \{\langle u, v \rangle \in R \, 
      | \, u \in L_H\}$ 
      and 
      $R \triangleright L_H = 
      \{\langle u, v \rangle \in R \, 
      | \, v \in L_H \}$. 
     We denote by $L_{H_1}$ and $L_{H_2}$
     the  languages $L_{H_1} = \psi^{-1} (H_1)
     \subset L_H$ and 
     $L_{H_2} = \psi^{-1} (H_2) \subset L_H$.             
     \begin{thm}       
     \label{nonautth1}        
        Assume that there exist 
        some FA--recognizable relations 
        $R_0$, $R_1$, $R_2$ $\subseteq S^* 
        \times S^*$  
        for which        
        $L_{H} \triangleleft R_0 
         = R_\varphi$,       
        $R_1 \triangleright L_{H_1}       
         = R_{p_1}$,
        $L_H \triangleleft R_2 
         = R_{p_2}$ and 
         $R_2 \triangleright \{w_0\} = 
         R_{p_2} \triangleright \{w_0\}$. 
        Then, for the function 
        $h(n)= 
         \max\{ d_A (\pi (w), \psi (w)) | 
         w \in L^{\leqslant n} \}$, 
         $\mathfrak{e} \preceq h$.          
         In particular, 
         $h \npreceq f$ for any $f \prec 
         \mathfrak{e}$.      
     \end{thm}  
     \beginpf
        Let $\eta (a,b,c)$ be the following 
        first--order formula: 
           \begin{equation*}
           \begin{split}            
            \eta (a,b,c)
            \equiv 
            \exists r,s_1,s_2,           
            t_1, t_2, t_3        
            \{ R_1(r,a) \land              
            (
                R_0 (b,s_1) \land 
                R_2 (s_1,s_2) \land                
                R_2(r,s_2)                                    
            )          
            \land  \\
            ( R_0 (r,t_1)  \land 
                   R_2 (t_1,t_2)) \land 
                  (R_2(c,w_0) \land 
                   R_0 (c,t_3) \land                   
                   R_2 (t_3,t_2)) \}.                     
            \end{split}
           \end{equation*}   
     Let us verify that for 
     any $a,b \in L_{H_1}$ 
     the formula $\eta(a,b,c)$ is true 
     if and only if $c \in L_{H_1}$ and 
     $\psi(a)+ \psi(b) = \psi(c)$ in the 
     cyclic group $H_1$.         
      Suppose that, 
      for some $a,b \in L_{H_1}$,      
      $\eta(a,b,c)$ 
      is true.  Let
           $\psi (a) =\left[\begin{array}{c}
       k \\     
       0 
       \end{array}
      \right] $, 
      $\psi (b) =\left[\begin{array}{c}
       \ell \\     
       0 
       \end{array}
      \right]$ for some $k,\ell \in \mathbb{Z}$.       
      Since $R_1(r,a)$ is true and 
      $R_1 \triangleright L_{H_1} = R_{p_1}$,
      then  $r \in L_H$ and 
      $\psi(r) = \left[\begin{array}{c}
       k \\     
       \star 
       \end{array}
      \right]$.
      Furthermore,  since $R_0 (b,s_1) \land 
                R_2 (s_1,s_2) \land                
                R_2(r,s_2)$ is true
       and $L_{H} \triangleleft R_0 = R_\varphi$,
       $L_H \triangleleft R_2 = 
       R_{p_2}$,  
       then $s_1 \in L_H$, $s_2 \in L_{H_2}$
       and 
       $\psi(s_1) = \left[\begin{array}{c}
       \ell \\     
       \ell
       \end{array}
      \right]$, $\psi(s_2) = 
      \left[\begin{array}{c}
       0 \\     
       \ell
       \end{array}
      \right]$,
       $\psi(r)= \left[\begin{array}{c}
       \star \\     
       \ell
       \end{array}
      \right]$. 
      Therefore, 
      $\psi(r)= \left[\begin{array}{c}
        k \\     
        \ell
       \end{array}
      \right]$. 
     Moreover, since 
     $R_0(r,t_1) \land R_2(t_1,t_2)$ 
     is true and $L_H \triangleleft R_0 = 
     R_\varphi$, $L_H \triangleleft R_2 = 
     R_{p_2}$, then   
     $\psi (t_1) = \left[\begin{array}{c}
        k \\     
       k+ \ell
       \end{array}
      \right]$,  $\psi(t_2) =
      \left[\begin{array}{c}
        0 \\     
       k+ \ell
       \end{array}
      \right] 
      $. 
      Finally, since 
      $R_2(c,w_0) \land R_0 (c,t_3) 
       \land R_2(t_3,t_2)$ and 
      $R_2 \triangleright \{w_0\} = 
       R_{p_2} \triangleright \{w_0\}$,
      $L_H \triangleleft R_0 = R_\varphi$, 
      $L_H \triangleleft R_2 = R_{p_2}$, then 
      $c=\left[\begin{array}{c}
        m \\     
        0 
       \end{array}
      \right]$, 
      $\psi(t_3) =
      \left[\begin{array}{c}
        m \\     
        m 
       \end{array}
      \right]$,
      $\psi(t_2) =
      \left[\begin{array}{c}
        0 \\     
        m 
       \end{array}
      \right]$.
      Thus, $c \in L_{H_1}$ and  
      $m = k + \ell$ which implies that
      $\psi (a) + \psi(b) = \psi (c)$.         
      The reverse is straightforward.


     Let $R \subseteq S^* \times S^* 
     \times S^*$ be the relation defined  
     by $\eta$, that is, $R(a,b,c)$ is true 
     iff $\eta(a,b,c)$ is true.       
     Since $R_0,R_1,R_2$ are FA--recognizable, 
     $R$ is FA--recognizable. 
     Let $(M,\times)$ be a monoid generated by $s$, 
     where $M = \{s^n \,|\, n \geqslant 0\} 
     \subset H_1$ and $\times$ is 
     the group multiplication in  $H_1$. 
     Clearly, 
     $(M,\times) \cong (\mathbb{N},+)$.    
     Let $u_n = \psi^{-1} (s^n) \in L$      
     and $L_M = \{ u_n | n \geqslant 0\}
     \subset L$. 
     It follows directly from  
     \cite[Lemma~6]{NiesThomas08} 
     (this lemma was originally proved 
     in \cite{BakhAndreFrankSasha} for 
     automatic monoids) that 
     there exist constants $C,N_0$ for which 
     $|u_n| \leqslant C \log n$
     for all $n \geqslant N_0$, where
     $|u_n|$ is the length of the string $u_n$. 
     It follows from the metric inequalities for 
     the Heisenberg 
     group $\mathcal{H}_3 (\mathbb{Z})$, 
     see, e.g., 
     \cite[Proposition~1.38]{Roebook}, that 
     there exist a constant $C_1>0$ for which
     $d_A (s^n) \geqslant C_1 n$ for all 
     $n \geqslant 0$.   
     We have: 
     $d_A (\pi(u_n),\psi(u_n)) \geqslant  
      d_A (\psi(u_n)) - d_A (\pi (u_n)) 
      \geqslant C_1 n - |u_n| \geqslant 
      C_1 n - C \log n $
     for all $n \geqslant N_0$. 
     Therefore, there exist some 
     constants $C_2>0$ and $N_1 \geqslant N_0$ 
     for which 
     $d_A (\pi(u_n),\psi(u_n)) \geqslant
      C_2 n \geqslant C_2 
      \exp \left( \frac{1}{C} |u_n|\right)$ 
      for all $n \geqslant N_1$.
     Clearly, the set
     $\{|u_n| \, | \, u_n \in L_M\} \subseteq
     \mathbb{N}$ is 
     infinite. Moreover, 
     by the finite difference lemma 
     (see, e.g., \cite[Lemma~14.1]{KKM11})  
     , $||u_{n+1}| -|u_n|| 
     \leqslant D$ for every $n \geqslant 0$ 
     and some constant $D$. Therefore, 
     there exists a constant $D_0$
     such that for every $j \geqslant D_0$ 
     there is $u_n \in L_M$ for which
     $j \geqslant |u_n| > j - D$. Thus, 
     for the function 
     $h(n)=\max\{d_A(\pi(w),\psi(w))|
     w \in L^{\leqslant n}\}$, 
     $\mathfrak{e} \preceq h$. 
     The last statement of the theorem 
     is straightforward.                     
   \endpf   
   \begin{rem}
      We note that the conditions 
      of Theorem \ref{nonautth1} 
      are clearly satisfied for the Cayley 
      automatic representation of 
      the Heisenberg group 
      $\mathcal{H}_3(\mathbb{Z})$ 
      constructed in 
      \cite[Section~6]{measuring_closeness1}. 
      As for FA--recognizable 
      relation $R_0 \subset S^* \times S^*$ 
      for which $L_H \triangleleft R_0 =  
                 R_\varphi$, it exists 
      if, for example, one
      additionally requires that 
      the left multiplication by $p^{-1}$
      in the group $\mathcal{H}_3 (\mathbb{Z})$      
      is FA--recognizable; it follows
      from the fact that for any 
      $h \in H$: $p^{-1} h p = \varphi (h)$.  
   \end{rem}
        
  \section{\bf Linear Upper Bounds for Almost All 
  Elements in Groups of Exponential Growth} 
  \label{linearupperalmostall}  
  
    In this section we show 
    that for an arbitrary  
    bijection $\psi: L \rightarrow G$ 
    between a language 
    $L \subseteq (A \cup A^{-1})^*$ and 
    a group $G$ of exponential growth      
    a linear upper bound  
    $d_A \left(\pi\left(\psi^{-1}(g)\right),
    g\right) \leqslant C |\psi^{-1}(g)|$
    holds for almost all $g \in G$
    in a certain sense,    
    see Theorem 
    \ref{allmostall_quasigeodesic_theorem1} 
    and Remark 
    \ref{allmostall_quasigeodesic_remark1}.
    However, in the following Remark  
    \ref{exptowergrowthremark}
    we show how to construct Cayley automatic 
    representations of the lamplighter
    group $\mathbb{Z}_2 \wr \mathbb{Z}$ 
    for which the function 
    \eqref{def_func_h_1} grows  
    faster than any tower of exponents.  
 \begin{thm}
 \label{allmostall_quasigeodesic_theorem1}        
  Let us assume that  
  $\psi: L \rightarrow G$ 
  is a Cayley automatic representation of a group 
  $G$ which has exponential growth.      
  Then there exist constants 
  $\lambda_1,\lambda_2>0$ 
  such that for almost all $g \in G$:  
  $\lambda_1 d_A (g) \leqslant 
   |w| \leqslant \lambda_2 d_A(g)$, 
  where $\psi(w) = g$. The term 
  almost all here means that 
  $\lim\limits_{n \rightarrow \infty}
   \frac{\#Q_n}{\# B_n} = 1$, where 
  $B_n = \{g\in G \, | \,
    d_A (g) \leqslant n \}$ is the ball of 
  radius $n$ in $G$  
  and  $Q_n \subseteq B_n$ is defined 
   as $Q_n = \{g \in B_n \,|\,  
        \lambda_1 d_A (g) \leqslant 
       |w| \leqslant \lambda_2 d_A(g)\}$.     
  In particular, for every $g \in  Q_n$, 
       $d_A(\pi(w),\psi(w)) \leqslant 
       \left( 1 + 
       \frac{1}{\lambda_1}\right)|w|$.  
 \end{thm}
 \beginpf
    The inequality 
    $|w| \leqslant \lambda_2 d_A (g)$ 
    always holds for some $\lambda_2>0$ 
    due to the bounded difference lemma.     
    Since $G$ has exponential growth, 
    there exists $\lambda>1$ for which 
    $\# B_n \geqslant \lambda^n$ for 
    all $n \geqslant n_0$. 
    For a given integer 
    $k>0$ we denote by $R_k$ the following 
    finite subset of $G$: 
    $R_k = \{g \in G \,|\,g=\psi(w), 
     w \in L^{< k} \}$; where 
     $L^{<k}=\{w \in L\,|\, |w|<k\}$. 
    Since $L \subseteq S^*$, 
    $\# R_k \leqslant |S|^{k-1}$.
    We denote by 
    $T_{n,k}$ the set 
    $T_{n,k} = B_n \setminus R_k$.
    For every $g \in T_{n,k}$, 
    $|w| \geqslant k$. Therefore, 
    if $\lambda_1 \leqslant \frac{k}{n}$,  
    for every $g \in T_{n,k}$ we have
    that $\lambda_1 d_A(g) \leqslant |w|$;  
    so $T_{n,k} \subseteq Q_n$.      
    We notice that 
    $\frac{\#T_{n,k}}{\# B_n} \geqslant
    1 - \frac{\#R_k}{\# B_n} 
    \geqslant 1- \frac{|S|^{k-1}}{\lambda^n}$ 
    for all $n \geqslant n_0$.
    So, it is enough to provide $\lambda_1$ 
    and a sequence 
    $k_n, n\geqslant n_0$ for which 
    $\lambda_1 \leqslant \frac{k_n}{n}$ for 
    all $n \geqslant n_0$ and 
    $\lim\limits_{n \rightarrow \infty}
     \frac{|S|^{k_n-1}}{\lambda^n}=0$. 
    We note that 
    $\frac{|S|^{k_n-1}}{\lambda^n}=
     \frac{|S|^{k_n-1}}{|S|^{ (\log_{|S|} 
     \lambda)n}} = \frac{1}{|S|^{(\log_{|S|} 
     \lambda)n - k_n + 1}}$. 
    Let us put $k_n = 
    \lceil \frac{1}{2}(\log_{|S|}\lambda)n
    \rceil$ for all $n \geqslant n_0$ 
    and $\lambda_1= 
    \frac{1}{2} (\log_{|S|}\lambda)$. Therefore, 
    $(\log_{|S|} \lambda)n - k_n +1 \geqslant 
    \frac{1}{2} (\log_{|S|}\lambda) n$, 
    so $\lim\limits_{n \rightarrow \infty}
     \frac{|S|^{k_n-1}}{\lambda^n} = 
     0$. Moreover, 
     $\frac{k_n}{n} \geqslant \lambda_1$
     for all $n\geqslant n_0$.           
     In order to prove the last inequality,
     we observe that 
     $d_A (\pi(w),\psi(w)) \leqslant 
      |w| + d_A (g) \leqslant |w| + 
      \frac{1}{\lambda_1}|w|$.  
 \endpf
 
 \begin{rem}
 \label{allmostall_quasigeodesic_remark1} 
    It is easy to see that Theorem 
    \ref{allmostall_quasigeodesic_theorem1} 
    holds for any bijection 
    $\psi:L \rightarrow G$ such that
    $||\psi^{-1}(ga)|-|\psi^{-1}(g)||
    \leqslant C$ for all $g \in G$ and 
    every generator $a \in A$, where  
    $C$ is a constant. 
    Moreover, for the inequality  
    $\lambda_1 d_A (g) \leqslant |w|$ and, 
    accordingly, the inequality 
    $d_A (\pi(w),\psi(w))\leqslant 
     \left(1+\frac{1}{\lambda_1}\right)|w|$, 
    no assumption is needed -- it 
    holds for almost all $g \in G$ 
    for any bijection $\psi$ between a language 
    $L$ and a group $G$ of exponential 
    growth. 
    Since in this paper 
    we focus mainly on Cayley automatic 
    representations of groups, in Theorem
    \ref{allmostall_quasigeodesic_theorem1}
    we assume that $\psi : L \rightarrow G$  
    is a Cayley automatic representation of $G$.                
 \end{rem}
 \begin{rem}
  \label{exptowergrowthremark}  
    We note that although for any 
    Cayley automatic representation 
    $\psi: L \rightarrow G$ of 
    a group of exponential growth $G$ 
    the inequality $d_A(\pi(w),\psi(w))
    \leqslant C |w|$      
    holds for some constant $C$    
    for almost all 
    $\psi(w)=g \in G$     
    in the sense of
    Theorem  
    \ref{allmostall_quasigeodesic_theorem1},
    it does not hold for all $g \in G$.             
    For example, let us 
    consider the following 
    Cayley  automatic representation 
    $\varphi: L_\Sigma \rightarrow 
              \mathbb{Z}_2 \wr \mathbb{Z}$ 
    over the alphabet 
    $\Sigma = \{+,-,0,1,C_0,C_1,\#\}$. 
    For any given pair $(f,z) \in 
                    \mathbb{Z}_2 \wr 
                    \mathbb{Z}$, we 
                    represent it as the 
    string: $u \# f(s) \dots  C_{f(z)} 
               \dots  f(t)$, where 
            $s$ and $t$ are the minimum
            and the maximum integers of 
            the set 
            $\{i| f(i)=1\} \cup \{z\}$, 
            $C_{f(z)}$ is $C_0$ or $C_1$ 
            if $f(z)=0$ or $f(z)=1$, 
            respectively, and the string $u$ is 
            a binary representation of
            the integer $s$. 
    For example, let us consider a 
    pair $(f,3)$, where $f(1)=1$,$f(2)=1$ and 
    $f(i) = 0$ if $i \neq 1,2$, it is 
    represented as the string:
    $+1\#11C_0$. Let us consider a pair 
    $(f,-3)$, where $f(-4)=1$,$f(-3)=1$,
    $f(-2)=1$, $f(i)=0$ if $i \neq -4,-3,-2$,
    it is represented as the string:
     $-100\#1C_11$. 
    We also refer the reader to 
    \cite[Example~4.2.1]{BerdinskyCSThesis}.
    
    One can then convert this 
    representation $\varphi$, 
    in a same way as in Lemma \ref{facts_lemma1},
    into some representation 
    $\psi: L \rightarrow 
     \mathbb{Z}_2 \wr \mathbb{Z}$ over the 
    alphabet $S=A \cup A^{-1}$, 
    where $A = \{a,t \}$ is the 
    standard set of generators of 
    $\mathbb{Z}_2 \wr \mathbb{Z}$: 
    $a$ is the nontrivial element 
    of $\mathbb{Z}_2$ and $t$ is a generator
    of $\mathbb{Z}$ (here we treat 
    $\mathbb{Z}_2$ and $\mathbb{Z}$
    as the subgroups of 
    $\mathbb{Z}_2 \wr \mathbb{Z}$).  
    Although 
    $\mathbb{Z}_2 \wr \mathbb{Z} 
     \in \mathcal{B}_\mathfrak{i}$ 
     \cite[Theorem~13]{measuring_closeness1}, 
    for the  representation 
    $\psi$ the inequality 
    $d_A (\pi(w),\psi(w)) \leqslant 
     C |w|$ does not hold for all
    $w \in L$ and any constant $C$. 
    In order to see that, 
    let us consider the representatives  
    $w_i =\psi^{-1}(g_i)$    
    of the elements 
    $g_i = (f_0,i) \in 
     \mathbb{Z}_2 \wr \mathbb{Z}$, $i>0$ 
     with respect to $\psi$, where 
    $f_0 (j) = 0$ for all $j \in \mathbb{Z}$. 
    Apparently,
    $d_A (g_i)=i$
    but the function 
    $\ell(i) = |w_i|$    
    grows, 
    coarsely, 
    as $\log i$. 
    So, the function 
    $h(n)=\mathrm{max}
     \{ d_A (\pi(w),\psi (w))| 
     w \in L ^{\leqslant n}\}$ grows 
    at least as fast as the exponential 
    function.           
    
    Moreover  one can construct a Cayley 
    automatic representation 
    $\psi: L \rightarrow 
     \mathbb{Z}_2 \wr \mathbb{Z}$ 
    for    which the function 
    $h (n) = \mathrm{max}
     \{ d_A (\pi(w),\psi (w))| 
     w \in L^{\leqslant n}\}$ 
     grows faster than any tower of exponents 
     $e^{e^{\dots ^e}}$. This follows 
     from the result shown by Frank Stephan:
    \end{rem}     
     \begin{thm}[Frank Stephan
          \cite{Stephan_email1}]
        There exists 
        an automatic representation 
        $\tau: 
        L_\tau \rightarrow \mathbb{N}$ 
        of the structure $(\mathbb{N},S)$,
        where $S$ is the successor function, 
        for which the function 
        $r(n) = \max \{\tau (w)| 
        w \in L_{\tau}^{\leqslant n} \}$
        grows faster than any  
        tower of exponents 
        $e^{e^{\dots ^e}}$. 
     \end{thm}              
    Clearly, one cannot directly 
    generalize Theorem 
    \ref{allmostall_quasigeodesic_theorem1}
    for Cayley automatic groups of subexponential 
    growth. Moreover, it simply does not
    hold for many Cayley automatic groups of subexponential
    growth -- consider, for example, a
    binary representation of the infinite cyclic
    group $\mathbb{Z}$.
     A f.g. group of 
    subexponential growth has either 
    intermediate growth or polynomial growth.      
    Miasnikov and Savchuk 
    constructed a FA--presentable graph 
    of intermediate growth  
    \cite{MiasnikovSavchuk13}. 
    However, it is still 
    unknown whether there exists 
    any Cayley automatic group of intermediate
    growth. As for f.g. 
    groups of polynomial growth,
    due to celebrated 
    Gromov's theorem \cite{GromovPolGrowth}, 
    any such group is virtually nilpotent. 


\vspace{0.3cm}\noindent{\bf Acknowledgment}

The authors thank 
  Murray Elder, Bakhadyr Khoussainov and 
  Frank Stephan for useful comments.

\label{lastpage}

\end{document}